\documentclass[a4paper,11pt,UKenglish]{article}
\usepackage[top=3.2cm, bottom=3cm, left=3.45cm, right=3.45cm]{geometry}
\usepackage[T1]{fontenc}
\usepackage[utf8]{inputenc}
\usepackage{latexsym}
\usepackage{amsmath,amsthm,amssymb,amsfonts}
\usepackage{mathtools}
\usepackage{scalerel}
\usepackage[english]{babel}
\usepackage{pstricks}
\usepackage{rotating}
\usepackage{graphicx}

\usepackage{tikz}
\usetikzlibrary{patterns, matrix}
\tikzstyle{NE-lines}=[pattern=north east lines, pattern color=black!45]

\newtheorem{theorem}{Theorem}[section]

\newtheorem{lemma}[theorem]{Lemma}
\newtheorem{corollary}[theorem]{Corollary}
\newtheorem{openpr}{Open Problem}
\newtheorem{remark}{Remark}

\newcommand{\cvd}{\hfill $\blacksquare$\bigskip}
\newcommand{\U}{\mathtt{U}} 
\newcommand{\D}{\mathtt{D}} 
\newcommand{\out}{s_{\sigma}}
\newcommand{\Cay}{\mathcal{C}}
\newcommand{\mapsigma}{\mathcal{S}^{\sigma}}
\newcommand{\pathsigma}{\mathcal{P}_{\sigma}}

\newcommand{\Sort}{\mathrm{Sort}}
\newcommand{\HPS}{\mathcal{H}\mathcal{P}\mathcal{S}}
\newcommand{\TPS}{\mathcal{T}\mathcal{P}\mathcal{S}}

\def\sdwys #1{\xHyphenate#1$\wholeString}
\def\xHyphenate#1#2\wholeString {\if#1$%
\else\say{\ensuremath{#1}}\hspace{2pt}%
\takeTheRest#2\ofTheString
\fi}
\def\takeTheRest#1\ofTheString\fi
{\fi \xHyphenate#1\wholeString}
\def\say#1{\begin{turn}{-90}\ensuremath{#1}\end{turn}}

\newenvironment{onestack}
{
	\begin{footnotesize}
	\psset{xunit=0.0355in, yunit=0.0355in, linewidth=0.02in}
	\begin{pspicture}(-15,-5)(40,30)
	\psline{c-c}(-15,25)(7.5,25)(7.5,0)(17.5,0)(17.5,25)(40,25)
	\rput[l](-15.5,22.5){\mbox{output}}
	\rput[r](40,22.5){\mbox{input}}
}
{
	\end{pspicture}
	\end{footnotesize}
}

\newcommand{\fillonestack}[4]{%
	\rput[l](-15,27.5){\ensuremath{#1}}
	\rput[c](12.6, 12.5){\begin{sideways}{\sdwys{#2}}\end{sideways}}
	\rput[r](40,27.5){\ensuremath{#3}}
	\rput[h](-10,0){\textbf{Step} \ensuremath{\boldsymbol{#4}}}
}

\begin{document}

\title{Sorting Cayley Permutations \\ with Pattern-avoiding Machines \thanks{The author is member of the INdAM Research group GNCS; he is partially supported by INdAM - GNCS 2019 project ``Propriet\'a combinatorie e rilevamento di pattern in strutture discrete lineari e bidimensionali".}}

\author{Giulio Cerbai \\ Dipartimento di Matematica e Informatica "U. Dini" \\ University of Firenze, Italy \\ giulio.cerbai@unifi.it}


\maketitle

\begin{abstract}
Pattern avoiding machines were recently introduced by Claesson, Ferrari and the current author to gain a better understanding of the classical $2$-stacksort problem. In this paper we generalize these devices by allowing permutations with repeated elements, also known as Cayley permutations. The main result is a description of those patterns such that the corresponding set of sortable permutations is a class. We also show a new involution on the set of Cayley permutations, obtained by regarding a pattern-avoiding stack as an operator. Finally, we analyze two generalizations of pop-stack sorting on Cayley permutations. In both cases we describe sortable permutations in terms of pattern avoidance.
\end{abstract}

\section{Introduction}\label{section_intro}

The problem of sorting a permutation using a stack, together with its many variants, has been widely studied in the literature. The original version was proposed by Knuth in~\cite{Kn}: given an input permutation $\pi$, either \textit{push} the next element of $\pi$ into the stack or \textit{pop} the top element of the stack, placing it into the output. The goal is to describe and enumerate sortable permutations. To sort a permutation means to produce a sorted output, i.e. the identity permutation. An elegant answer can be given in terms of pattern avoidance: a permutation is sortable if and only if it does not contain a subsequence of three elements which is order isomorphic to $231$. A set of permutations that can be characterized in terms of pattern avoidance is called a \textit{class} and the minimal excluded permutations are its \textit{basis}. The notion of pattern avoidance turns out to be a fundamental tool to approach a great variety of problems in combinatorics. We refer the reader to~\cite{B} for a more detailed survey on stack-sorting disciplines, and to~\cite{B2} and~\cite{Ki} for an overview on patterns in permutations and words. It is easy to realize that the optimal algorithm for the classical stacksorting problem has two key properties. First, the elements in the stack are maintained in increasing order, reading from top to bottom. Moreover, the algorithm is \emph{right-greedy}, meaning that it always performs a push operation, unless this violates the previous condition. Note that the expression "right-greedy" refers to the usual (and most natural) representation of this problem, depicted in Figure~\ref{stacksort_machine}.

Although the classical problem is rather simple, as soon as one allows several stacks connected in series things become much harder. For example, it is known that the permutations that can be sorted using two stacks in series form a class, but in this case the basis is infinite~\cite{M}, and still unknown. The enumeration of such permutations is still unknown too. In the attempt of gaining a better understanding of this device, some (simpler) variants have been considered. A \textit{pop-stack} is a (restricted) stack where all the elements are extracted everytime a pop operation is performed. Pop-stacks were introduced by Avis and Newborn~\cite{AN}, where the authors prove that permutations sortable through a pop-stack are the so called layered permutations. A permutation is \textit{layered} if it avoids $231$ and $312$. More recently, two or more pop-stacks in series were considered in~\cite{AS} and \cite{CG}. In his PhD thesis~\cite{W}, West considered two passes through a classical stack, which is equivalent to perform a \textit{right-greedy} algorithm on two stacks in series. In~\cite{Sm}, Smith considered a decreasing stack followed by an increasing stack. This machine was then generalized in~\cite{CCiF} to the case of many decreasing stacks, followed by an increasing one. Recently, the authors of~\cite{CCF} considered \textit{pattern-avoiding machines}, an even more general device consisting of two stacks in series with a right-greedy procedure, where a restriction on the first stack is given in terms of pattern avoidance. More precisely, the first stack is not allowed to contain an occurrence of a forbidden pattern $\sigma$, for a fixed $\sigma$. West's device is obtained by choosing $\sigma=21$. The pattern $\sigma=12$ corresponds to the device analyzed in~\cite{Sm}, but with a right-greedy (and thus less powerful) algorithm.

Other than imposing restrictions on devices and sorting algorithms, one can also allow a larger set of input sequences. Since the notion of pattern itself is inherently more general, it is natural to consider sorting procedures on bigger sets of strings~\cite{AAAHH,ALW,DK}. Here we pursue this line of research by analyzing the behaviour of pattern-avoiding machines on permutations with repeated letters, which are known as Cayley permutations. A more formal definition of Cayley permutation will be given in Section~\ref{section_tools}, together with the necessary background and tools.

In Section~\ref{section_pattern_machine} we generalize a result of~\cite{CCF} by determining for which patterns $\sigma$ the words that can be sorted by the $\sigma$-machine form a class. In such cases, we also give an explicit description of the basis, which is either a singleton or consists of two patterns.

In Section~\ref{section_operator}, we regard a $\sigma$-avoiding stack as a function $\mapsigma$ that maps an input word into the resulting output, characterizing the patterns $\sigma$ that give rise to a bijective operator. The proof of this result relies on the encoding of $\mapsigma$ as a labeled Dyck path. By composing $\mapsigma$ with the reverse operator, we obtain a new involution on the set of Cayley permutations. Such involution has the nice property of preserving the multiset of entries of a Cayley permutation. It also leads to a constructive description of the set of sortable permutations.

In Section~\ref{section_popstack} we analyze two generalizations of pop-stack sorting on Cayley permutations. We call them hare pop-stack and tortoise pop-stack, in analogy with a paper by Defant and Kravitz~\cite{DK}. In both cases, we characterize the set of sortable permutations in terms of pattern avoidance. A simple geometrical description allows us to enumerate the set of tortoise pop-stack sortable permutations, while the hare case is left for a future investigation.

\begin{center}
\begin{figure}
\hspace{1.5cm}
\begin{minipage}{5cm}
\begin{footnotesize}
\psset{xunit=0.0355in, yunit=0.0355in, linewidth=0.02in}
\begin{pspicture}(0,0)(35,20)
\psline{c-c}(-10,15)(7.5,15)(7.5,-5)(17.5,-5)(17.5,15)(35,15)
\rput[l](-10.5,12.5){\mbox{output}}
\rput[r](35,12.5){\mbox{input}}
\psline[linearc=0.2,arrows=->](22.5,20)(17,18.5)(14.5,13.5)
\psline[linearc=0.2,arrows=->](10.5,13.5)(8,18.5)(2.5,20)
\rput[t](19,23.5){\ensuremath{push}}
\rput[t](6,23.5){\ensuremath{pop}}
\rput[r](35,17.5){\ensuremath{\pi_1 \dots \pi_n}}
\rput[h](12.5,-10){\ensuremath{\left\lfloor \begin{array}{c}
                          2 \\
                          1
                        \end{array}
\right\rfloor}}	
\end{pspicture}
\end{footnotesize}
\end{minipage}
\hspace{1cm}
\begin{minipage}{5cm}
\begin{footnotesize}
\psset{xunit=0.0355in, yunit=0.0355in, linewidth=0.02in}
\begin{pspicture}(0,0)(35,20)
\psline{c-c}(-15,15)(2.5,15)(2.5,-5)(12.5,-5)(12.5,15)(22.5,15)(22.5,-5)(32.5,-5)(32.5,15)(50,15)
\rput[l](-15.5,12.5){\mbox{output}}
\rput[r](50,12.5){\mbox{input}}
\rput[h](27.5,-10){\ensuremath{\sigma}}
\rput[h](7.5,-10){\ensuremath{\left\lfloor \begin{array}{c}
                          2 \\
                          1
                        \end{array}
\right\rfloor}}
\rput[l](-25,17.5){\ensuremath{}}
\rput[c](12.6, 2.5){\begin{sideways}{\sdwys{}}\end{sideways}}
\rput[c](22.6, 2.5){\begin{sideways}{\sdwys{}}\end{sideways}}
\rput[r](50,17.5){\ensuremath{\pi_1 \dots \pi_n}}
\psline[linearc=0.2,arrows=->](37.5,20)(32,18.5)(29.5,13.5)
\psline[linearc=0.2,arrows=->](5.5,13.5)(3,18.5)(-2.5,20)
\psline[linearc=0.2,arrows=->](25.5,14)(20.5,18.8)(14.5,18.5)(9.5,14)
\rput[t](34,23.5){\ensuremath{push_{1}}}
\rput[t](1,23.5){\ensuremath{pop}}
\rput[t](17.5,23.5){\ensuremath{push_{2}}}
\end{pspicture}
\end{footnotesize}
\end{minipage}
\vspace{1.2cm}
\caption{Sorting with one stack (on the left) and sorting with two stacks, where the first one is $\sigma$-restricted (on the right).}
\label{stacksort_machine}
\end{figure}
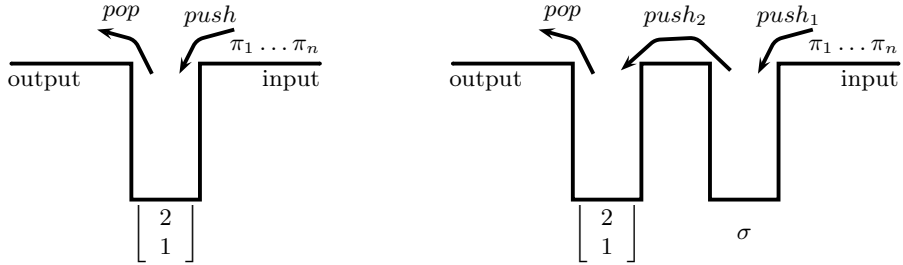
\end{center}

\section{Tools and Notations}\label{section_tools}

Let $\mathbb{N}^*$ be the set of strings over the alphabet $\mathbb{N}=\left\lbrace 1,2,\dots \right\rbrace$ of positive integers. Let $x=x_1 \cdots x_n$ and $p=p_1 \cdots p_k$ in $\mathbb{N}^*$, with $k \le n$. The word $x$ \textit{contains} the pattern $p$ if there are indices $i_1<i_2<\cdots <i_k$ such that $x_{i_1} x_{i_2} \cdots x_{i_k}$ is order isomorphic to $p$. Equivalently, for each pair of indices $u,v$, $x_{i_u}<x_{i_v}$ if and only if $p_u<p_v$ and $x_{i_u}=x_{i_v}$ if and only if $p_u=p_v$. In this case, we write $p \le x$ and we say that $x_{i_1} x_{i_2} \cdots x_{i_k}$ is an \textit{occurrence} of $p$. Otherwise, we say that $x$ \textit{avoids} the pattern $p$. This notion generalizes the usual notion of pattern involvement on permutations. For example, the string $x=142215$ contains the pattern $2113$, since the substring $4225$ is order isomorphic to $2113$. On the other hand, $x$ avoids the pattern $1234$. A \textit{class} is a set of words that is closed downwards with respect to pattern involvement. A class is determined by the minimal set of words it avoids, which is called its \textit{basis}.

Denote by $\Cay$ the set of strings $\pi$ on $\mathbb{N}$ where each integer from $1$ to $max(\pi)$ appears at least once. Following~\cite{MF}, we call these strings \textit{Cayley permutations} (they are called \textit{normalized words} in~\cite{DK}, and sometimes also surjective words, Fubini words or packed words). Cayley permutations, with respect to their length, are enumerated by sequence $A000670$ in the OEIS~\cite{Sl}. For example, the only Cayley permutation of length one is the string $1$, there are three Cayley permutations of length two, namely $11$, $12$ and $21$, and thirteen Cayley permutations of length three, which are $111$, $112$, $121$, $211$, $122$, $212$, $221$, $123$, $132$, $213$, $231$, $312$, $321$. Since only the relative order of the elements is relevant for avoidance and containment, patterns live naturally in the set $\Cay$. More precisely, given $x \in \mathbb{N}^*$, an order-isomorphic string $\pi \in \Cay$ can be produced by suitably rescaling the elements of $x$, so to remove gaps. For this reason, and because we think that the most natural setting is the one where patterns and words belong to the same set, in the rest of the paper we will work on $\Cay$ rather than on $\mathbb{N}^*$. We denote by $\Cay(\pi)$ the set of Cayley permutations that avoid the pattern $\pi$, for $\pi \in \Cay$; for a set of patterns $B= \{ \pi_1,\dots,\pi_k \}$, the set of Cayley permutations that avoid every pattern $\pi_1,\cdots,\pi_k$ is denoted by $\Cay(B)$. The \textit{reverse} of the Cayley permutation $\pi=\pi_1\cdots\pi_k$ is $\pi^r=\pi_k\cdots\pi_1$. A \textit{weak descent} is a pair of consecutive elements $\pi_i,\pi_{i+1}$ such that $\pi_i \ge \pi_{i+1}$. If $\pi > \pi_{i+1}$, the pair is a \textit{(strong) descent}. \textit{Weak ascents} and \textit{(strong) ascents} are defined analogously.

\section{$\sigma$-machines on Cayley Permutations}\label{section_pattern_machine}

The authors of~\cite{CCF} introduced pattern-avoiding machines on permutations. Here we generalize these devices by allowing Cayley permutations both as inputs and as forbidden patterns. Let $\sigma$ be a Cayley permutation of length at least two. A $\sigma$\textit{-stack} is a stack that is not allowed to contain an occurrence of the pattern $\sigma$ when reading the elements from top to bottom. Before introducing $\sigma$-machines, we recall some useful results. Classical stacksort on $\mathbb{N}^*$ has been discussed in~\cite{DK}. Note that there are two possibilities when defining the analogue of the stacksort algorithm on $\mathbb{N}^*$. One can either allow a letter to sit on a copy of itself in the stack, or force a pop operation if the next element of the input is equal to the top element of the stack. Here we choose the former possibility, leaving the latter for future investigation. This is equivalent to regard a classical stack as a $21$-avoiding stack. The following theorem, proved in~\cite{DK} for $\mathbb{N}^*$, also applies to Cayley permutations. 

\begin{theorem}\label{hare_stacksort}
Let $\pi$ be a Cayley permutation. Then $\pi$ is sortable using a $21$-stack if and only if $\pi$ avoids $231$.
\end{theorem}

The term $\sigma$\textit{-machine} refers to performing a right-greedy algorithm on two stacks in series: a $\sigma$-stack, followed by a $21$-avoiding stack (see Figure~\ref{stacksort_machine}). A Cayley permutation $\pi$ is $\sigma$\textit{-sortable} if the output of the $\sigma$-machine on input $\pi$ is weakly increasing. The set of $\sigma$-sortable permutations is denoted by $\Sort(\sigma)$. We use the notation $\out(\pi)$ to denote the output of the $\sigma$-stack on input $\pi$. Note that, since $\out(\pi)$ is the input of the $21$-stack, Theorem~\ref{hare_stacksort} guarantees that $\pi \in \Sort(\sigma)$ if and only if $\out(\pi)$ avoids $231$. This fact will be used repeatedly for the rest of the paper. In~\cite{CCF}, the authors provided a characterization of the (permutation) patterns $\sigma$ such that the set of $\sigma$-sortable permutations is a class. The main goal of this section is to extend this result to Cayley permutations.

\begin{remark}\label{remark_hat}
Let $\sigma=\sigma_1 \cdots \sigma_k$ be a Cayley permutation. If an input Cayley permutation $\pi$ avoids $\sigma^r$, then the restriction of the $\sigma$-stack is never triggered and $\out(\pi)=\pi^r$. Otherwise, the leftmost occurrence of $\sigma$ results necessarily in an occurrence of $\hat{\sigma}$ in $\out(\pi)$, where $\hat{\sigma}=\sigma_2 \sigma_1 \sigma_3 \sigma_4 \cdots \sigma_k$. 
\end{remark}

From now on, we denote by $\hat{\sigma}$ the Cayley permutation obtained from $\sigma$ by interchanging $\sigma_1$ and $\sigma_2$.

\begin{theorem}\label{suff_class}
Let $\sigma$ be a Cayley permutation. If $\hat{\sigma}$ contains $231$, then $\Sort(\sigma)=\Cay(132,\sigma^r)$. In this case, $\Sort(\sigma)$ is a class with basis either $\{ 132,\sigma^r \}$, if $\sigma^r$ avoids $132$, or $\{ 132\}$, otherwise.
\end{theorem}
\proof We start by proving that $\Sort (\sigma )\subseteq \Cay (132,\sigma^r )$. Let $\pi \in \Sort (\sigma)$. Equivalently, suppose that $\out(\pi)$ avoids $231$. Suppose by contradiction that $\pi$ contains $\sigma^r$. Then $\out(\pi)$ contains $\hat{\sigma}$ due to Remark~\ref{remark_hat} and $\hat{\sigma}$ contains $231$ by hypothesis, which is impossible. Otherwise, if $\pi$ avoids $\sigma^r$, but contains $132$, then $\out(\pi)=\pi^r$ due to the same remark. Moreover $\pi^r$ contains $231$ by hypothesis, again a contradiction with $\pi \in \Sort(\sigma)$. This proves that $\Sort (\sigma )\subseteq \Cay (132,\sigma^r )$.\\
Conversely, suppose that $\pi$ avoids both $132$ and $\sigma^r$. Then $\out(\pi)=\pi^r$, which avoids $132^r=231$ by hypothesis, therefore $\pi$ is $\sigma$-sortable. This completes the proof.
\cvd

Next we show that the condition of Theorem~\ref{suff_class} is also necessary for $\Sort(\sigma)$ in order to be a class. The only exception is given by the pattern $\sigma=12$.

\begin{theorem}\label{12_stack}
$\Sort(12)=\Cay(213)$.
\end{theorem}
\proof Let $\pi$ be a Cayley permutation. Suppose that the element $1$ appears $k$ times in $\pi$ and write $\pi =A_1 1 A_2 1 \cdots A_k 1 A_{k+1}$.  It is easy to see that
$$s_{12}(\pi )=s_{12}(A_1) s_{12}(A_2) \cdots s_{12}(A_k) s_{12} (A_{k+1}) 1 \cdots 1.$$
Indeed a copy of $1$ can enter the $12$-stack only if the $12$-stack is either empty or it contains only other copies of $1$. Finally, the element $1$ cannot play the role of $2$ in an occurrence of the (forbidden) pattern $12$. Therefore the presence of some copies of $1$ at the bottom of the $12$-stack does not affect the sorting process of the block $A_i$, for each $i$.\\
Now, suppose that $\pi$ contains an occurrence $bac$ of $213$. We prove that $\pi$ is not $12$-sortable by showing that $s_{12}(\pi)$ contains $231$. We proceed by induction on the length of $\pi$. Write $\pi =A_1 1 A_2 1 \cdots A_k 1 A_{k+1}$ as above. Suppose that $b \in A_i$ and $c \in A_j$, for some $i \le j$ (note that $b,c \neq 1$). If $i=j$, then $A_i$ contains an occurrence $bac$ of $213$. Thus $s_{12}(A_i)$ contains $231$ by induction, as wanted. Otherwise, let $i<j$. Then $b \in s_{12}(A_i)$ and $c \in s_{12}(A_j)$ and the elements $b$ and $c$, together with any copy of $1$, realize an occurrence of $231$ in $s_{12}(\pi)$, as desired.\\
Conversely, suppose that $\pi=\pi_1 \cdots \pi_n$ is not sortable, i.e. $s_{12}(\pi )$ contains $231$. We prove that $\pi$ contains $213$. Let $bca$ an occurrence of $231$ in $s_{12}(\pi)$. Note that $b$ has to precede $c$ in $\pi$. This is due to the fact that a non-inversion in the output necessarily comes from a non-inversion in the input, since the stack is $12$-avoiding. However, $b$ is pushed out before $c$ enters. Denote with $x$ the next element of the input when $b$ is extracted. Then we have $x<b$ and also $x \neq c$, since $c>b$. Finally, the triple $bxc$ forms an occurrence of $213$ in $\pi$, as desired.
\cvd

\begin{theorem}\label{necess_class}
Let $\sigma$ be a Cayley permutation and suppose $\sigma \neq 12$. If $\hat{\sigma}$ avoids $231$, then $\Sort (\sigma)$ is not a class.
\end{theorem}
\proof Let $\sigma=\sigma_1 \cdots \sigma_k$, with $k \ge 2$. We show that there are two Cayley permutations $\alpha,\beta$ such that $\alpha$ contains $\beta$, $\beta$ is $\sigma$-sortable and $\alpha$ is not $\sigma$-sortable. This proves that $\Sort(\sigma)$ is not closed downwards, as desired. Figure~\ref{table_alpha_beta} shows an example of such $\alpha$ and $\beta$ for patterns $\sigma$ of length two and for $\sigma=231$. Now, suppose that $\sigma$ has length at least three and $\sigma \neq 231$. Then the Cayley permutation $\alpha=132$ is not $\sigma$-sortable. Indeed, $\out(\alpha)=\alpha^r=231$, since $\alpha$ avoids $\sigma^r$. Next we define the permutation $\beta$ according to the following case by case analysis.

\begin{figure}
\begin{center}
\begin{tabular}{|c|c|c|}
  \hline
  $\sigma$ & $\alpha \notin \Sort(\sigma)$ & $\beta \ge \alpha, \beta \in \Sort(\sigma)$ \\
  \hline
  11 & 132 & 3132 \\
  \hline
  21 & 132 & 35241 \\
  \hline
  231 & 1324 & 361425 \\
  \hline
\end{tabular}
\end{center}
\caption{The case by case analysis of Theorem~\ref{necess_class}.}\label{table_alpha_beta}
\end{figure}

\begin{figure}
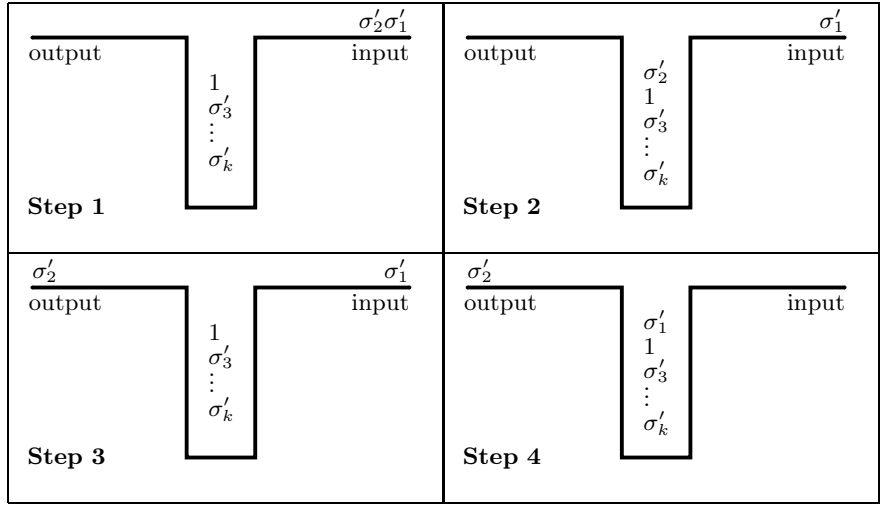

\begin{center}
\begin{tabular}{|c|c|}
\hline
\begin{onestack}
\fillonestack{}{{\sigma'_k}{.}{.}{.}{\sigma'_3}{1}}{\sigma'_2 \sigma'_1}{1}
\end{onestack}
&
\begin{onestack}
\fillonestack{}{{\sigma'_k}{.}{.}{.}{\sigma'_3}{1}{\sigma'_2}}{\sigma'_1}{2}
\end{onestack}
\\
\hline
\begin{onestack}
\fillonestack{{\sigma'_2}}{{\sigma'_k}{.}{.}{.}{\sigma'_3}{1}}{\sigma'_1}{3}
\end{onestack}
&
\begin{onestack}
\fillonestack{{\sigma'_2}}{{\sigma'_k}{.}{.}{.}{\sigma'_3}{1}{\sigma'_1}}{}{4}
\end{onestack}
\\
\hline
\end{tabular}
\end{center}
\caption{The action of the $\sigma$-stack on input $\beta$ described in the proof of Theorem~\ref{necess_class}.}\label{sorting_alpha}
\end{figure}

\begin{itemize}
\item Suppose that $\sigma_1$ is the strict minimum of $\sigma$, i.e. $\sigma_1=1$ and $\sigma_i \ge 2$ for each $i \ge 2$. Define
$$\beta=\sigma'_k \cdots \sigma'_3 1 \sigma'_2 \sigma'_1,$$
where $\sigma'_i=\sigma_i+1$ for each $i$. Note that $\beta$ is a Cayley permutation and $1 \sigma'_2 \sigma'_1$ is an occurrence of $132$ in $\beta$. We prove that $\beta$ is $\sigma$-sortable by showing that $\out(\beta)$ avoids $231$. The action of the $\sigma$-stack on input $\beta$ is depicted in Figure~\ref{sorting_alpha}. The first $k-1$ elements of $\beta$ are pushed into the $\sigma$-stack, since $\sigma$ has length $k$. Then the $\sigma$-stack contains $1 \sigma'_3 \cdots \sigma'_k$, reading from top to bottom, and the next element of the input is $\sigma'_2$. Note that $\sigma'_2>1$, whereas $\sigma_1<\sigma_2$, therefore $\sigma'_2 1 \sigma'_3 \cdots \sigma'_k$ is not an occurrence of $\sigma$ and $\sigma'_2$ is pushed. The next element of the input is now $\sigma'_1$. Here $\sigma'_1 \sigma'_2 \sigma'_3 \cdots \sigma'_k$ is an occurrence of $\sigma$, thus we have to pop $\sigma'_2$ before pushing $\sigma'_1$. After the pop operation, the $\sigma$-stack contains $1 \sigma'_3 \cdots \sigma'_k$. Again $\sigma'_1>1$, whereas $\sigma_1<\sigma_2$, therefore $\sigma'_1$ is pushed. The resulting string is
$$\out(\beta)=\sigma'_2 \sigma'_1 1 \sigma'_3 \sigma'_4 \cdots \sigma'_k.$$
We wish to show that $\out(\beta)$ avoids $231$. Note that $\sigma'_2 \sigma'_1 \sigma'_3 \sigma'_4 \cdots \sigma'_k \simeq \hat{\sigma}$ avoids $231$ by hypothesis. Moreover, the element $1$ cannot be part of an occurrence of $231$, because $\sigma'_2 > \sigma'_1$ and $1$ is strictly less than the other elements of $\beta$. Therefore $\out(\beta)$ avoids $231$, as desired.

\item Otherwise, suppose that $\sigma_1$ is not the strict minimum of $\sigma$, i.e. either $\sigma_1 \neq 1$ or $\sigma_i=1$ for some $i \ge 2$. Define
$$\beta=\sigma''_k \cdots \sigma''_2 1 \sigma''_1 2,$$
where $\sigma''_i=\sigma_i+2$ for each $i$. Note that $\beta$ is a Cayley permutation and $1 \sigma''_2 2$ is an occurrence of $132$ in $\beta$. Consider the action of the $\sigma$-stack on $\beta$. Again the first $k-1$ elements of $\beta$ are pushed into the $\sigma$-stack. Then the $\sigma$-stack contains $\sigma''_2 \cdots \sigma''_k$, reading from top to bottom, and the next element of the input is $1$. Note that $1 \sigma''_2 \cdots \sigma''_k$ is not an occurrence of $\sigma$. Indeed $1<\sigma''_i$ for each $i$, while $\sigma_1$ is not the strict minimum of $\sigma$ by hypothesis. Therefore $1$ enters the $\sigma$-stack. The next element of the input is then $\sigma''_1$, which realizes an occurrence of $\sigma$ together with $\sigma''_2 \cdots \sigma''_k$. Thus $1$ and $\sigma''_2$ are extracted before $\sigma''_1$ is pushed. Finally, the last element of the input is $2$. Again $2$ can be pushed into the $\sigma$-stack because $2$ is strictly smaller than every element in the $\sigma$-stack, whereas $\sigma_1$ is not the strict minimum of $\sigma$ by hypothesis. The resulting string is
$$\out(\beta) = 1 \sigma''_2 2 \sigma''_1 \sigma''_3 \cdots \sigma''_k.$$
Note that $\sigma''_2 \sigma''_1 \sigma''_3 \cdots \sigma''_k \simeq \hat{\sigma}$ avoids $231$ by hypothesis. Finally, it is easy to realize that the elements $1$ and $2$ cannot be part of an occurrence of $231$, similarly to the previous case. This completes the proof.
\end{itemize}
\cvd

\begin{corollary}\label{class_vs_nonclass}
Let $\sigma$ be a Cayley permutation of length three or more. Then the set of $\sigma$-sortable permutation $\Sort(\sigma)$ is not a class if and only if $\hat{\sigma}$ avoids $231$. Otherwise, if $\hat{\sigma}$ contains $231$, then $\Sort(\sigma)$ is a class with basis either $\{ 132,\sigma^r \}$, if $\sigma^r$ avoids $132$, or $\{ 132\}$, otherwise.
\end{corollary}

We end this section by analyzing the $21$-machine. The $11$-machine will be discussed in Section~\ref{section_operator}, thus completing the analysis of the $\sigma$-machines on Cayley permutations for patterns $\sigma$ of length two. The classical permutation analogue of the $21$-machine is exactly the (well known) case of the West $2$-stack sortable permutations~\cite{W2}. In this case, although sortable permutations do not form a class, it is possible to describe them efficiently in terms of avoidance of barred patterns.

\begin{theorem}\cite{W2}\label{west_2stack}
A permutation $\pi$ is not $21$-sortable if and only if $\pi$ contains $2341$ or $\pi$ contains an occurrence of the barred pattern $3 \bar{5}241$, i.e. an occurrence $3241$ which is not part of an occurrence of $35241$.
\end{theorem}

\begin{figure}
\begin{center}
$\mathcal{W} \;=\;
\begin{tikzpicture}[scale=0.6, baseline=20pt]
\fill[NE-lines] (1,4) rectangle (2,5);
\draw [semithick] (0,1) -- (5,1);
\draw [semithick] (0,2) -- (5,2);
\draw [semithick] (0,3) -- (5,3);
\draw [semithick] (0,4) -- (5,4);
\draw [semithick] (1,0) -- (1,5);
\draw [semithick] (2,0) -- (2,5);
\draw [semithick] (3,0) -- (3,5);
\draw [semithick] (4,0) -- (4,5);
\filldraw (1,3) circle (6pt);
\filldraw (2,2) circle (6pt);
\filldraw (3,4) circle (6pt);
\filldraw (4,1) circle (6pt);
\end{tikzpicture}
$
\qquad
$\mathcal{Z} \;=\;
\begin{tikzpicture}[scale=0.6, baseline=20pt]
\fill[NE-lines] (1.15,3.85) rectangle (1.85,4.15);
\fill[NE-lines] (1.15,4.15) rectangle (1.85,5);
\draw [semithick] (0,0.85) -- (5,0.85);
\draw [semithick] (0,1.15) -- (5,1.15);
\draw [semithick] (0,1.85) -- (5,1.85);
\draw [semithick] (0,2.15) -- (5,2.15);
\draw [semithick] (0,2.85) -- (5,2.85);
\draw [semithick] (0,3.15) -- (5,3.15);
\draw [semithick] (0,3.85) -- (5,3.85);
\draw [semithick] (0,4.15) -- (5,4.15);
\draw [semithick] (0.85,0) -- (0.85,5);
\draw [semithick] (1.15,0) -- (1.15,5);
\draw [semithick] (1.85,0) -- (1.85,5);
\draw [semithick] (2.15,0) -- (2.15,5);
\draw [semithick] (2.85,0) -- (2.85,5);
\draw [semithick] (3.15,0) -- (3.15,5);
\draw [semithick] (3.85,0) -- (3.85,5);
\draw [semithick] (4.15,0) -- (4.15,5);
\filldraw (1,3) circle (6pt);
\filldraw (2,2) circle (6pt);
\filldraw (3,4) circle (6pt);
\filldraw (4,1) circle (6pt);
\end{tikzpicture}
$
\end{center}
\caption{On the left, the barred pattern $3 \bar{5}241$, which is equivalent to the mesh pattern $\mathcal{W}= (3241,\{(1,4) \})$. The shaded box keeps into account the case of an occurrence of $3241$ that is part of a $35241$. On the right, the Cayley-mesh pattern $\mathcal{Z}$. The additional shaded region in $\mathcal{Z}$ keeps into account the case of an occurrence of $3241$ that is part of an occurrence of $34241$.}\label{mesh_west}
\end{figure}
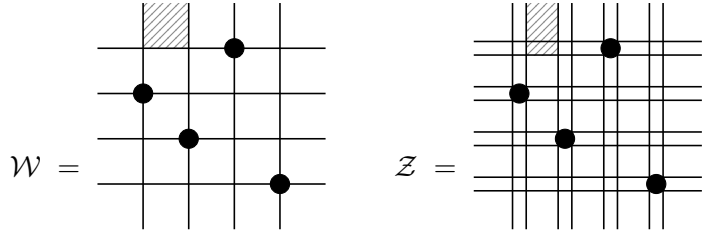

The previous theorem can be reformulated in terms of a more general notion of pattern, which will be useful later when dealing with Cayley permutations. A \textit{mesh pattern}~\cite{BC} of length $k$ is a pair $(\tau,A)$, where $\tau$ is a permutation of length $k$ and $A \subseteq \left[ 0,k \right] \times \left[ 0,k \right]$ is a set of pairs of integers. The elements of $A$ identify the lower left corners of shaded squares in the plot of $\tau$ (see Figure~\ref{mesh_west}). An occurrence of the mesh pattern $(\tau,A)$ in the permutation $\pi$ is then an occurrence of the classical pattern $\tau$ in $\pi$ such that no other elements of $\pi$ are placed into a shaded square of $A$.

Note that the barred pattern $3 \bar{5} 241$ is equivalent to the mesh pattern $\mathcal{W}$ depicted in Figure~\ref{mesh_west}. Now, in order to prove an analogous characterization for the $12$-machine on Cayley permutations, we need to adapt the definition of mesh pattern to strings that may contain repeated elements. In other words, we allow the shading of regions that correspond to repeated elements. Instead of giving a formal definition, we refer to the example illustrated in Figure~\ref{mesh_west}. We will use the term \textit{Cayley-mesh pattern} to denote mesh patterns on Cayley permutations.

\begin{lemma}\label{increas_stack}
Let $\pi=\pi_1 \cdots \pi_n$ be a Cayley permutation. Suppose that $\pi_i < \pi_j$, for some $i < j$. Then $\pi_i$ precedes $\pi_j$ in $s_{21}(\pi)$.
\end{lemma}
\proof It follows from the definition of $21$-stack.
\cvd

\begin{theorem}
A Cayley permutation $\pi$ is not $21$-sortable if and only if $\pi$ contains $2341$ or $\pi$ contains the Cayley-mesh pattern $\mathcal{Z}$ depicted in Figure~\ref{mesh_west}. In particular, $\Sort(21)$ is not a class. For example, the $21$-sortable Cayley permutation $34241$ contains the non-sortable pattern $3241$.
\end{theorem}
\proof We can basically repeat the argument used by West for classical permutations. The only difference is the additional shaded box, which corresponds to an occurrence of $3241$ that is part of an occurrence of $34241$. We sketch the proof anyway for completeness.\\
Let $\pi$ be a Cayley permutation and suppose that $\pi$ is $21$-sortable. Suppose by contradiction that $\pi$ contains an occurrence $bcda$ of $2341$ and consider the action of the $21$-stack on $\pi$. By Lemma~\ref{increas_stack}, $b$ is extracted from the $21$-stack before $c$ enters. Similarly, $c$ is extracted before $d$ enters. Thus $s_{21}(\pi)$ contains the occurrence $bca$ of $231$, against $\pi$ sortable. Otherwise, suppose that $\pi$ contains an occurrence $cbda$ of $3241$. We show that there is an element $x$ between $c$ and $b$ in $\pi$ such that $x \ge d$. If $x<c$ for each $x$ in between $c$ and $b$, then $b$ is pushed into the $21$-stack before $c$ is popped. This results in the occurrence $bca$ of $231$ in $s_{21}(\pi)$, a contradiction with $\pi$ $21$-sortable. Otherwise, suppose there is at least one element $x$ between $c$ and $b$ in $\pi$, with $x \ge c$. If $x=c$, we can repeat the same argument with $xbda$ instead of $cbda$. If $c<x<d$, then $cxda \simeq 2341$, which is impossible due to what said in the above case. Therefore it has to be $x \ge d$, as desired.\\
Conversely, suppose that $\pi$ is not $12$-sortable. Equivalently, let $bca$ be an occurrence of $231$ in $s_{21}(\pi)$. We show that either $\pi$ contains $2341$ or $\pi$ contains an occurrence $cbda$ of $3241$ such that $x<d$ for each $x$ between $c$ and $b$ in $\pi$. Observe that $a$ follows $c$ and $b$ in $\pi$ due to Lemma~\ref{increas_stack}. Suppose that $b$ comes before $c$ in $\pi$. Note that $c$ is extracted from the $21$-stack before $a$ enters. Let $d$ the next element of the input when $c$ is extracted. Then $d>c$ and $bcda$ is an occurrence of $2341$, as wanted. Otherwise, suppose that $b$ follows $c$ in $\pi$, and thus $\pi$ contains $cba$. Since $c$ is not extracted before $b$ enters, it has to be $x \le c$ for each $x$ between $c$ and $b$ in $\pi$. Moreover, $c$ is extracted before $a$ enters. When $c$ is extracted, the next element $d$ of the input is such that $d>c$. This results in an occurrence $cbda$ of $3241$ with the desired propriety.
\cvd

\begin{openpr}
Enumerate the $21$-sortable Cayley permutations. The initial terms of the sequence are $1,3,13,73,483,3547,27939,231395$ (not in~\cite{Sl}).
\end{openpr}

\section{$\sigma$-stacks as Operators}\label{section_operator}

In this section we regard $\sigma$-stacks as operators. Let $\sigma$ be a Cayley permutation and define the map $\mapsigma: \Cay \mapsto \Cay$ by $\mapsigma(\pi)=\out(\pi)$, for each $\pi \in \Cay$. We are interested in the behavior of the map $\mapsigma$. This line of inquiry for stacksort operators is not new in the literature. More generally, suppose to perform a deterministic sorting procedure. Then it is natural to consider the map $\mathcal{S}$ that associates to an input string $\pi$ the (uniquely determined) output of the sorting process. Some of the arising problems are the following.
\begin{itemize}
\item Determine the \textit{fertility} of a string, which is the number of its pre-images under $\mathcal{S}$. Fertilty under classical stacksort has been recently investigated by Defant~\cite{D}.
\item Determine the image of $\mathcal{S}$, i.e. the strings with positive fertility. These are often called \textit{sorted permutations} ~\cite{BM}.
\end{itemize}

We start by discussing the case $\sigma=11$. Here we provide a useful decomposition that allows us to determine explicitly the image $\mapsigma(\pi)$ of any input Cayley permutation $\pi$. From now on, we denote by $\mathcal{R}$ the \textit{reverse} operator, i.e. $\mathcal{R}(\pi)=\pi^r$, for each $\pi \in \Cay$. 

\begin{lemma}\label{11_decom}
Let $\sigma=11$ and let $\pi=\pi_1 \cdots \pi_n$ be a Cayley permutation. Suppose that $\pi$ contains $k+1$ occurrences $\pi_1$,$\pi_1^{(1)}$,$\dots$,$\pi_1^{(k)}$ of $\pi_1$, for some $k \ge 0$. Write $\pi=\pi_1 B_1 \pi_1^{(1)} B_2 \cdots \pi_1^{(k)} B_k$. Then
$$\mathcal{S}^{11}(\pi)= \mathcal{S}^{11}(B_1) \pi_1 \mathcal{S}^{11}(B_2) \pi_1^{(1)} \cdots \mathcal{S}^{11}(B_k) \pi_1^{(k)}.$$
\end{lemma}
\proof Consider the action of the $11$-stack on input $\pi$. Since $x \neq \sigma_1$ for each $x \in B_1$, the sorting process of $B_1$ is not affected by the presence of $\sigma_1$ at the bottom of the $11$-stack. Then, when the next element of the input is the second occurrence $\sigma_1^{(1)}$ of $\sigma_1$, the $11$-stack is emptied, since $\sigma_1 \sigma_1^{(1)}$ is an occurrence of the forbidden $11$. The first elements of $\mathcal{S}^{11}(\pi)$ are thus $\mathcal{S}^{11}(B_1) \sigma_1$. Finally, $\sigma_1^{(1)}$ is pushed into the (empty) $11$-stack and the same argument can be repeated.
\cvd

\begin{theorem}\label{11_bij}
Let $\sigma=11$. Then $(\mathcal{R} \circ \mathcal{S}^{11})$ is an involution on $\Cay$. Moreover, $\mathcal{S}^{11}$ is a length-preserving bijection on $\Cay$. Therefore, the number of $11$-sortable Cayley permutations of length $n$ is equal to the number of $231$-avoiding Cayley permutations of length $n$.
\end{theorem}
\proof We proceed by induction on the length of the input permutation. Let $\pi=\pi_1 \cdots \pi_n$ a Cayley permutation of length $n$. The case $n=1$ is trivial. If $n \ge 2$, write $\pi=\pi_1 B_1 \pi_1^{(1)} B_2 \cdots \pi_1^{(k)} B_k$ as in the previous lemma. Then, using the same lemma and the inductive hypothesis:
\begin{equation*}
\begin{split}
\left[ \mathcal{R} \circ \mathcal{S}^{11} \right]^2(\pi)=& \\
\left[ \mathcal{R} \circ \mathcal{S}^{11} \right]^2 \left( \pi_1 B_1 \pi_1^{(1)} B_2 \cdots \pi_1^{(k)} B_k \right) =&\\
\left[ \mathcal{R} \circ \mathcal{S}^{11} \circ \mathcal{R} \right] \left( \mathcal{S}^{11}(B_1) \pi_1 \mathcal{S}^{11}(B_2) \pi_1^{(1)} \cdots \mathcal{S}^{11}(B_k) \pi_1^{(k)} \right) =& \\
\left[ \mathcal{R} \circ \mathcal{S}^{11} \right] \left( \pi_1^{(k)} \mathcal{R}(\mathcal{S}^{11}(B_k)) \cdots \pi_1^{(1)} \mathcal{R}(\mathcal{S}^{11}(B_2)) \pi_1 \mathcal{R}(\mathcal{S}^{11}(B_1)) \right) =& \\
\mathcal{R} \left(  \mathcal{S}^{11}(\mathcal{R}(\mathcal{S}^{11}(B_k))) \pi_1^{(k)} \cdots \mathcal{S}^{11}(\mathcal{R}(\mathcal{S}^{11}(B_2))) \pi_1^{(1)} \mathcal{S}^{11}(\mathcal{R}(\mathcal{S}^{11}(B_1))) \pi_1 \right) =& \\
\pi_1 \left[ \mathcal{R} \circ \mathcal{S}^{11} \right]^2 (B_1) \pi_1^{(1)} \left[ \mathcal{R} \circ\mathcal{S}^{11} \right]^2(B_2) \cdots \pi_1^{(k)} \left[ \mathcal{R} \circ \mathcal{S}^{11} \right]^2(B_k) =& \\
\pi_1 B_1 \pi_1^{(1)} B_2 \cdots \pi_1^{(k)} B_k = \pi & \\
\end{split}
\end{equation*}
Therefore $(\mathcal{R} \circ \mathcal{S}^{11})^2(\pi)=\pi$, as desired. Finally, the reverse map $\mathcal{R}$ is bijective, thus $\mathcal{S}^{11}$ is a bijection on $\Cay$ with inverse $\mathcal{R} \circ \mathcal{S}^{11} \circ \mathcal{R}$.
\cvd

Theorem~\ref{11_bij} provides a constructive description of the set $\Sort(11)$. Indeed, since $\Sort(11)=\mathcal{R} \circ \mathcal{S}^{11} \circ \mathcal{R} (\mathcal{C}(231))$, every $11$-sortable permutation $\pi$ is obtained from a $231$-avoiding Cayley permutation by applying $\mathcal{R} \circ \mathcal{S}^{11} \circ \mathcal{R}$. Next we generalize the above result by providing a characterization of all patterns $\sigma$ such that $\mapsigma$ is bijective on $\Cay$. The main tool is an encoding of the action of $\mapsigma$ as a Dyck path.

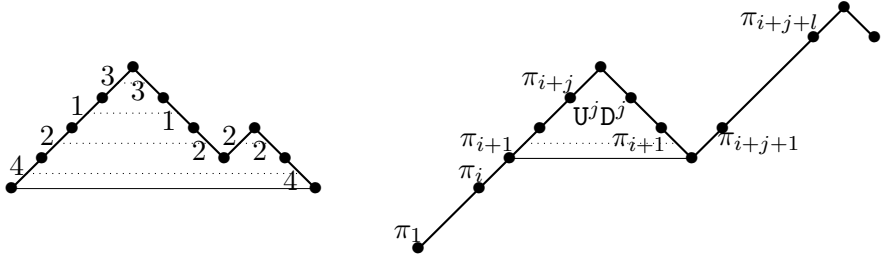
\begin{figure}
\begin{minipage}{5cm}
\begin{tikzpicture}[scale=0.4]
\draw [ultra thin] (0,0) -- (10,0);
\draw [thick] (0,0) -- (4,4);
\draw [thick] (4,4) -- (7,1);
\draw [thick] (7,1) -- (8,2);
\draw [thick] (8,2) -- (10,0);
\draw [dotted] (0.5,0.5) -- (9.5,0.5);
\draw [dotted] (1.5,1.5) -- (6.5,1.5);
\draw [dotted] (2.5,2.5) -- (5.5,2.5);
\draw [dotted] (3.5,3.5) -- (4.5,3.5);
\draw [dotted] (7.5,1.5) -- (8.5,1.5);
\node[] at (0,0) {$\bullet$};
\node[] at (1,1) {$\bullet$};
\node[] at (2,2) {$\bullet$};
\node[] at (3,3) {$\bullet$};
\node[] at (4,4) {$\bullet$};
\node[] at (5,3) {$\bullet$};
\node[] at (6,2) {$\bullet$};
\node[] at (7,1) {$\bullet$};
\node[] at (8,2) {$\bullet$};
\node[] at (9,1) {$\bullet$};
\node[] at (10,0) {$\bullet$};
\node[above,left] at (0.75,0.75) {$4$};
\node[above,left] at (1.75,1.75) {$2$};
\node[above,left] at (2.75,2.75) {$1$};
\node[above,left] at (3.75,3.75) {$3$};
\node[below,left] at (4.75,3.25) {$3$};
\node[below,left] at (5.75,2.25) {$1$};
\node[below,left] at (6.75,1.25) {$2$};
\node[above,left] at (7.75,1.75) {$2$};
\node[below,left] at (8.75,1.25) {$2$};
\node[below,left] at (9.75,0.25) {$4$};
\end{tikzpicture}
\end{minipage}
\begin{minipage}{5cm}
\begin{tikzpicture}[scale=0.4]
\draw [thick] (0,0) -- (6,6);
\draw [thick] (6,6) -- (9,3);
\draw [thick] (9,3) -- (14,8);
\draw [thick] (14,8) -- (15,7);
\node[] at (0,0) {$\bullet$};
\node[] at (2,2) {$\bullet$};
\node[] at (3,3) {$\bullet$};
\node[] at (4,4) {$\bullet$};
\node[] at (5,5) {$\bullet$};
\node[] at (6,6) {$\bullet$};
\node[] at (7,5) {$\bullet$};
\node[] at (8,4) {$\bullet$};
\node[] at (9,3) {$\bullet$};
\node[] at (10,4) {$\bullet$};
\node[] at (13,7) {$\bullet$};
\node[] at (14,8) {$\bullet$};
\node[] at (15,7) {$\bullet$};
\node[] at (6,4.5) {$\U^j \D^j$};
\node[above,left] at (2.5,2.5){$\pi_{i}$};
\node[above,left] at (5.5,5.5){$\pi_{i+j}$};
\node[above,left] at (13.5,7.5){$\pi_{i+j+l}$};
\node[above,left] at (0.5,0.5){$\pi_{1}$};
\node[above,left] at (3.5,3.5){$\pi_{i+1}$};
\node[left] at (8.5,3.5){$\pi_{i+1}$};
\node[right] at (9.5,3.5){$\pi_{i+j+1}$};
\draw [dotted] (3.5,3.5) -- (8.5,3.5);
\draw [dotted] (1,1) -- (2,2);
\draw [thin] (3,3) -- (9,3);
\end{tikzpicture}
\end{minipage}
\caption{On the left, the Dyck path $UUUUDDDUDD$ which encodes $\mathcal{S}^{11}(42132)$. On the right, the (prefix of the) path $P_{\sigma}(\pi)$ mentioned in the proof of Corollary~\ref{sigma_bij}. Dotted lines connect matching steps, which have the same label.}\label{figure_bij_cases}\label{Dyck_path}
\end{figure}

A \emph{Dyck path} is a path in the discrete plane $\mathbb{Z}\times \mathbb{Z}$ starting at the origin, ending on the $x$-axis, never falling below the $x$-axis and using two kinds of steps (of length $1$), namely up steps $\U=(+1,+1)$ and down steps $\D=(+1,-1)$. The \textit{height} of a step is its final ordinate. For each up step $\U$, there is a unique \textit{matching} step $\D$ defined as the first $\D$ step after $\U$ with height one less than $\U$. The \emph{length} of a Dyck path is the total number of its steps. A \textit{valley} of a Dyck path is an occurrence of two consecutive steps $\D \U$. An example of Dyck path is illustrated in Figure~\ref{Dyck_path}. It is well known that Dyck paths, according to the semilength, are enumerated by Catalan numbers (sequence $A000108$ in~\cite{Sl}). A \textit{labeled Dyck path} is a Dyck path where each step has a label. In this paper we consider labeled Dyck paths where the label of each up step is the same as the label of its matching down step. Therefore we can represent a labeled Dyck path $\mathcal{P}$ as a pair $\mathcal{P}=(P,\pi)$, where $P$ is the underlying Dyck path and $\pi$ is the string obtained by reading the labels of the up steps from left to right. Given an unlabeled Dyck path $P$ of length $2n$, the \textit{reverse} path $\mathcal{R}(P)$ of $P$ is obtained by taking the symmetric path with respect to the vertical line $x=n$.

Now let $\sigma$ be a Cayley permutation and suppose we are applying $\mapsigma$ to the input Cayley permutation $\pi$, i.e. we are sorting $\pi$ using a $\sigma$-stack. Then define a labeled Dyck path $\pathsigma(\pi)$ as follows.
\begin{itemize}
\item Insert an up step $\U$ labeled $a$ whenever the algorithm pushes an element $a$ into the $\sigma$-stack.
\item Insert a down step $\D$ labeled $a$ whenever the algorithm pops an element $a$ from the $\sigma$-stack.
\end{itemize}

Equivalently, define $P_{\sigma}(\pi)$ as the unlabeled Dyck path obtained by recording the push operations of the $\sigma$-stack with $\U$ and the pop operations with $\D$. Then $\pathsigma(\pi)=(P_{\sigma}(\pi),\pi)$. Note that $P_{\sigma}(\pi)$ is a Dyck path. Indeed the number of push and pop operations performed when processing $\pi$ is the same, therefore the number of $\U$ steps matches the number of $\D$ steps (and thus the path ends on the $x$-axis). Moreover, the path cannot go below the $x$-axis, since this would correspond to performing a pop operation when the $\sigma$-stack is empty, which is not possible. An example of this construction, when $\sigma=11$, is depicted in Figure~\ref{Dyck_path}. Some basic properties of $\pathsigma(\pi)$ are listed in the following Lemma, whose straightforward proof is omitted.

\begin{lemma}\label{path_prop}
Let $\sigma$ be a Cayley permutation. Let $\pi=\pi_1 \cdots \pi_n$ be a Cayley permutation of length $n$ and let $\pathsigma(\pi)=(P_{\sigma}(\pi),\pi)$. Then:
\begin{enumerate}
\item The input $\pi$ is obtained by reading the labels of the up steps of $P_{\sigma}(\pi)$ from left to right. The output $\out(\pi)$ is obtained by reading the labels of the down steps from left to right.

\item The height of each up (respectively down) step of $P_{\sigma}(\pi)$ is equal to the number of elements contained in the $\sigma$-stack after having performed the corresponding push (respectively pop) operation.

\item The $\sigma$-stack is emptied by a pop operation if and only if the corresponding $\D$ step of $P_{\sigma}(\pi)$ is a return on the $x$-axis. In other words, the decomposition of $\pi$ considered in Lemma~\ref{11_decom} corresponds to the decomposition of $P_{\sigma}(\pi)$ obtained by considering the returns on the $x$-axis.

\item The labels of the down steps are uniquely determined by the labels of the up steps. Conversely, the labels of the down steps uniquely determine the labels of the up steps. More precisely, matching steps have the same label. Indeed the element pushed into the $\sigma$-stack by an up step is then popped by the matching down step.

\item Let $\D \U$ be a valley in $P_{\sigma}(\pi)$. Let $a$ be the label of $\D$ and $b$ the label of $\U$. Then $b$ plays the role of $\sigma_1$ in an occurrence of $\sigma$ that triggers the restriction of the $\sigma$-stack, whereas $a$ plays the role of $\sigma_2$ in such occurrence. Moreover the number of valleys of $P_{\sigma}(\pi)$ is equal to the number of elements of $\pi$ that trigger the restriction of the $\sigma$-stack.

\item If $\sigma_1=\sigma_2$, then, for each valley $\D \U$, the labels of $\D$ and $\U$ are the same.
\end{enumerate}
\end{lemma}

\begin{theorem}\label{dyck_reverse}
Let $\sigma=\sigma_1 \cdots \sigma_k$ be a Cayley permutation. Let $\pi=\pi_1 \cdots \pi_n$ be a Cayley permutation and let $\gamma=\mathcal{R}(\mapsigma(\pi))$. Consider the two labeled Dyck paths $\pathsigma(\pi)=(P_{\sigma}(\pi),\pi)$ and $\pathsigma(\gamma)=(P_{\sigma}(\gamma),\gamma)$.
\begin{enumerate}
\item If $\sigma_1=\sigma_2$, then $P_{\sigma}(\pi)=\mathcal{R} (P_{\sigma}(\gamma))$.
\item If $P_{\sigma}(\pi)=\mathcal{R} (P_{\sigma}(\gamma))$, then $(\mathcal{R} \circ \mapsigma )^2 (\pi)=\pi$.
\end{enumerate}
\end{theorem}
\proof 
\begin{enumerate}
\item Suppose that $\sigma_1=\sigma_2$. We proceed by induction on the number of valleys of $P_{\sigma}(\pi)$. If $P_{\sigma}(\pi)$ has zero valleys, then $\pi$ avoids $\mathcal{R}(\sigma)$ by item $5.$ of Lemma~\ref{path_prop}. Therefore $\mapsigma(\pi)=\mathcal{R}(\pi)$ and $\gamma=\mathcal{R}^2(\pi) = \pi$. Since $P_{\sigma}(\pi)=\U^n \D^n$ is a pyramid, and each pyramid is equal to its reverse, the thesis follows immediately.\\
Now suppose that $P_{\sigma}(\pi)$ has at least one valley. Let $P_{\sigma}(\pi)=p_1 \cdots p_{2n}$ and write $P_{\sigma}(\pi)=\U^i \U^j \D^j \U^l \D Q$, where the steps $p_{i+2j}$ and $p_{i+2j+1}$ form the leftmost valley and $Q=p_{i+2j+l+2} \cdots p_n$ is the remaining suffix of $P_{\sigma}(\pi)$ (see Figure~\ref{figure_bij_cases}). Note that the label of both $p_{i+2j}$ and $p_{i+2j+1}$ is equal to $\pi_{i+1}$ because of items $4.$, $5.$ and $6.$ of Lemma~\ref{path_prop}. Item $5.$ also implies that $p_{i+2j+1}$ plays the role of $\sigma_1$ in an occurrence of $\sigma$ that triggers the restriction of the $\sigma$-stack. More precisely, immediately after the push of $\pi_{i+j}$ (i.e. after the up step $p_{i+j}$ in $P_{\sigma}(\pi)$), $\pi_{i+j+1}$ is the next element of the input. Since the next segment of the path is $\D^j$, $j$ pop operations are performed before pushing $\pi_{i+j+1}$. This means that the element $\pi_{i+1}$, corresponding to the last down step, plays the role of $\sigma_2$ in an occurrence of $\sigma$, while $\pi_{i+j+1}$ plays the role of $\sigma_1$. Moreover there are $k-2$ elements in the $\sigma$-stack that play the role of $\sigma_3,\dots,\sigma_k$. Since the elements in the $\sigma$-stack correspond to the labels of the initial prefix $\U^i$, $\pi_1 \cdots \pi_i$ contains an occurrence of $\sigma_k \cdots \sigma_3$ (claim I). Then, after $j$ pop operations are performed, the $\sigma$-stack contains $\pi_i \cdots \pi_1$, reading from top to bottom, and the elements $\pi_{i+j+1},\pi_{i+j+2}, \dots, \pi_{i+j+l}$ are pushed (claim II). Now, write
$$\pi= \underbrace{\pi_1 \cdots \pi_i}_{A} \ \underbrace{\pi_i+1 \cdots \pi_{i+j}}_{B} \ \underbrace{\pi_{i+j+1} \cdots \pi_{i+j+l}}_{C} \ \underbrace{\pi_{i+j+l+1} \cdots \pi_n}_{D},$$
where the elements of $A$ correspond to the initial prefix $\U^i$ of $P_{\sigma}(\pi)$, $B$ corresponds to $\U^j$, $C$ to $\U^l$ and $D$ to the remaining up steps. Consider the string $A C D = \pi_1 \cdots \pi_i \pi_{i+j+1} \cdots \pi_n$ obtained by removing the segment $B=\pi_{i+1} \cdots \pi_{i+j}$ from $\pi$. Let $\tilde{\pi}$ be the only Cayley permutation that is order isomorphic to $A C D$ (i.e. obtained by suitably rescaling the elements of $A C D$, if necessary). Note that $\pathsigma(\tilde{\pi})$ is obtained from $\pathsigma(\pi)$ by cutting out the pyramid $\U^j \D^j$, which corresponds to the removed segment $B$. This is because the elements contained in the $\sigma$-stack after having pushed $\pi_i$ are exactly the same as the elements contained in the $\sigma$-stack after having pushed $\pi_{i+j+1}$, thus we can safely cut out the pyramid $\U^j \D^j$ without affecting the sorting procedure. Therefore
$$\mapsigma(\pi)=\mathcal{R}(B) \mapsigma(\tilde{\pi}) \quad \mbox{and} \quad \gamma=\mathcal{R}(\mapsigma(\pi))=\mathcal{R}(\mapsigma(\tilde{\pi}))B.$$
Now, since $P_{\sigma}(\tilde{\pi})$ has one valley less than $P_{\sigma}(\pi)$, by inductive hypothesis $P_{\sigma}(\tilde{\pi})=\mathcal{R}(P_{\sigma}( \tilde{\gamma}))$, where $\tilde{\gamma}=\mathcal{R}(\mapsigma(\tilde{\pi}))$. The only difference bewteen $P_{\sigma}(\pi)$ and $P_{\sigma}(\tilde{\pi})$ is the removed pyramid $\U^j \D^j$. If we show that $P_{\sigma}(\gamma)$ is obtained from $P_{\sigma}(\tilde{\gamma})$ by reinserting the same pyramid $\U^j \D^j$ in the same place, the thesis follows. We have $\gamma=\mathcal{R}(\mapsigma(\tilde{\pi}))B$ and $\tilde{\gamma}=\mathcal{R}(\mapsigma(\tilde{\pi}))$. Consider the last push performed by the $\sigma$-stack when processing $\tilde{\gamma}$, which corresponds to the last up step of $\pathsigma(\tilde{\gamma})$. Notice that, since $P_{\sigma}(\tilde{\pi})=\mathcal{R}(P_{\sigma}( \tilde{\gamma})$, this is also the first down step of $P_{\sigma}(\tilde{\pi})$, and thus the first pop performed when processing $\tilde{\pi}$. Therefore the elements contained in the $\sigma$-stack after the last push performed while processing $\tilde{\gamma}$ are $\pi_{i+j+l} \cdots \pi_{i+j+1} \pi_i \cdots \pi_1$, reading from top top bottom. If we sort $\gamma$ instead of $\tilde{\gamma}$, we have to process the additional segment $B$. Now, the first element of $B$ is $\pi_{i+1}$. Since the same happened when sorting $\pi$ (see claim I), $\pi_{i+1}$ realizes an occurrence of $\sigma$ together with $\pi_{i+j+1}$ (which plays the role of $\sigma_2$) and other $k-2$ elements in $\pi_1 \cdots \pi_i$. The only difference is that, contrary to what happened when sorting $\pi$, the role of $\pi_{i+1}$ and $\pi_{i+j+1}$ are interchanged: here the hypothesis $\sigma_1=\sigma_2$ is relevant. As a result, before pushing the first element $\pi_{i+1}$ of $B$, we have to pop each element of the $\sigma$-stack up to $\pi_{i+j+1}$, $\pi_{i+j+1}$ included. After that, the $\sigma$-stack contains $\pi_i \cdots \pi_1$, reading from top to bottom. Therefore we can push $\pi_{i+1}=\pi_{i+j+1}$ and the remaining elements of $B$ because of claim II. This means that $P_{\sigma}(\gamma)$ is obtained by inserting a pyramid $\U^j \D^j$ immediately before the last $i$ down steps of $P_{\sigma}(\tilde{\gamma})$, as desired.

\item By hypothesis, $P_{\sigma}(\gamma)=\mathcal{R}(P_{\sigma}(\pi))$, therefore the word $w$ obtained by reading the labels of the down steps of $P_{\sigma}(\gamma)$ (from left to right) is $w=\mathcal{R}(\pi)$. By definition of $\pathsigma(\gamma)$, we also have $w=\mapsigma(\gamma)$. Therefore $\mathcal{R}(\pi)=\mapsigma(\gamma)=\mapsigma(\mathcal{R}(\mapsigma(\pi)))$ and the thesis follows by applying the reverse operator to both sides of the equality.

\end{enumerate}
\cvd 

\begin{corollary}\label{sigma_bij}
Let $\sigma=\sigma_1 \cdots \sigma_k \in \Cay$. Then $\mapsigma$ is bijective if and only if $\sigma_1=\sigma_2$. In this case, $\mapsigma$ is a bijection on $\Cay$ that preserves the multiset of entries of a Cayley permutation and $\mathcal{R} \circ \mapsigma$ is an involution on the set $\Cay$.
\end{corollary}
\proof Suppose that $\sigma_1 \neq \sigma_2$. Then $\hat{\sigma} \neq \sigma$, thus also $\mathcal{R}(\sigma) \neq \mathcal{R}(\hat{\sigma})$. Finally, $\mapsigma(\mathcal{R}(\sigma))=\hat{\sigma}=\mapsigma((\mathcal{R}(\hat{\sigma})))$, therefore $\mapsigma$ is not injective.\\
Conversely, suppose that $\sigma_1=\sigma_2$. By Theorem~\ref{dyck_reverse}, we have that $(\mathcal{R} \circ \mapsigma)^2$ is the identity on $\Cay$, therefore $\mathcal{R} \circ \mathcal{\mapsigma}$ is bijective. Finally, since the reverse map $\mathcal{R}$ is bijective, $\mapsigma$ is a bijection too, as desired.
\cvd

\section{Pop-stack on Cayley permutations}\label{section_popstack}

This section is devoted to the study of pop-stack sorting on Cayley permutations. Recall from Section~\ref{section_intro} that a \textit{pop-stack} is a stack where all the elements are extracted everytime a pop operation is performed. In analogy with~\cite{DK}, we introduce the hare and tortoise variants of a pop-stack, according whether or not a letter is allowed to sit on a copy of itself.

A \textit{hare-popstack} is a $21$-popstack, i.e. a $21$-stack that is emptied everytime a pop operation is performed. A Cayley permutation $\pi$ is \textit{hare pop-stack sortable} if $\pi$ is sortable using a right-greedy algorithm on a hare-popstack.

A \textit{tortoise-popstack} is a $\{21,11 \}$-popstack, i.e. a $\{21,11 \}$-stack that is emptied everytime a pop operation is performed. A Cayley permutation $\pi$ is \textit{tortoise pop-stack sortable} if $\pi$ is sortable using a right-greedy algorithm on a tortoise-popstack.

Denote by $\HPS(\pi)$ and $\TPS(\pi)$ the output of a hare pop-stack and, respectively, a tortoise pop-stack, on input $\pi$. Recall that, since we are allowing repeated elements, to sort a Cayley permutation means to produce a weakly increasing Cayley permutation (and not necessarily the identity permutation as in the classical case). Equivalently, the output of either a hare or tortoise pop-stack is not sorted if and only if it contains a strong descent. We shall provide a characterization of hare and tortoise pop-stack sortable permutations in terms of forbidden patterns, starting with hare-popstack. The next lemma is a straightforward consequence of the definition of hare-popstack.

\begin{lemma}\label{blocks_dec_hare}
Let $\pi$ be a Cayley permutation. Write $\pi=B_1 B_2 \cdots B_k$, where each block $B_t$ is maximally weakly decreasing (i.e. the last element of each block $B_t$ forms a strong ascent together with first element of the next block $B_{t+1}$). Then
$$\HPS(x)=\mathcal{R} (B_1) \cdots \mathcal{R}(B_k).$$
\end{lemma}

\begin{theorem}\label{hare_pop_char}
Let $\pi$ be a Cayley permutation. Then $\pi$ is hare pop-stack sortable if and only if $\pi$ avoids $231$, $312$ and $2121$.
\end{theorem}
\proof Suppose that $\pi$ is hare pop-stack sortable. Observe that hare pop-stack sortable Cayley permutations are a subset of $12$-sortable Cayley permutations. Thus $\pi$ avoids $231$ due to Lemma~\ref{hare_stacksort}. Now suppose by contradiction that $\pi$ contains an occurrence $cab$ of $312$. Then, when $b$ enters the hare pop-stack, $a$ has been already extracted. Otherwise $\HPS(\pi)$ would not be weakly increasing, against the hypothesis. Therefore also $c$ has been extracted and thus $b<c$, which is again a contradiction with $\pi$ hare pop-stack sortable. Similarly, suppose that $\pi$ contains an occurrence $bab'a'$ of $2121$. Then $a$, and thus also $b$, must have been extracted before $b'$ enters the stack, since $b'>a$. Therefore $b$ is extracted before $a'$ enters the hare pop-stack, which is impossible because $a'<b$.\\
Conversely suppose that $\pi$ is not hare pop-stack sortable. We wish to show that $\pi$ contains an occurrence of either $231$, $312$ or $2121$. Write $\pi=B_1 B_2 \cdots B_k$ as in Lemma~\ref{blocks_dec_hare}. Then $\HPS(\pi)=B_1^R \cdots B_k^R$ and $\HPS(\pi)$ contains at least one strict descent. Let $a>b$ the leftmost strict descent in $\HPS(\pi)$. Due to Lemma~\ref{blocks_dec_hare}, it must be $a \in B_i$ and $b \in B_{i+1}$, for some $i$. The same result implies that $a$ is the first element of $B_i$ and $b$ is the last element of $B_{i+1}$. Now, denote by $u$ the last element of $B_i$ and by $v$ the first element of $B_{i+1}$, as illustrated in Figure~\ref{blocks_fig}. We have $a \ge u$, $v \ge b$ and $v>u$. Consider the following case by case analysis.

\begin{itemize}
\item If $B_i$ is a singleton, then $u=a>b$ and $a=u<v$, therefore $b \neq v$ and $avb$ is an occurrence of $231$ in $\pi$.

\item If $B_{i+1}$ is a singleton, then $b=v>u$ and $v=b<a$, therefore $a \neq u$ and $auv$ is an occurrence of $312$ in $\pi$.

\item Finally, suppose that both $B_i$ and $B_{i+1}$ are not singletons and consider the four elements $auvb$ in $\pi$. If $a>v$, then $auv$ is an occurrence of $312$. If $a<v$, then $avb$ is an occurrence of $231$. Otherwise, suppose that $a=v$ (and so $a=v>u$). Then $auvb$ is an occurrence of $2121$, if $u=b$; $aub$ is an occurrence of $231$, if $u>b$; and $aub$ is an occurrence of $312$, if $u<b$.
\end{itemize}
\cvd

\begin{figure}
\begin{equation*}
\begin{split}
\pi =& \underbrace{\cdots}_{B_1} | \cdots | \underbrace{a \cdots u}_{B_i} | \underbrace{v \cdots b}_{B_{i+1}} | \cdots | \underbrace{\cdots}_{B_k} \\
\mathcal{H} / \TPS(\pi) =& \underbrace{\cdots}_{B_1} | \cdots | \underbrace{u \cdots a}_{B_i} | \underbrace{b \cdots v}_{B_{i+1}} | \cdots | \underbrace{\cdots}_{B_k}
\end{split}
\end{equation*}
\caption{The decomposition of $\pi$ used in Theorem~\ref{hare_pop_char} and Theorem~\ref{tortoise_pop_char}.}\label{blocks_fig}
\end{figure}

Next we consider tortoise pop-stack.

\begin{lemma}\label{blocks_dec_tortoise}
Let $\pi$ be a Cayley permutation. Write $\pi=B_1 B_2 \cdots B_k$, where each block $B_t$ is maximally strictly decreasing (i.e. the last element of each block $B_t$ forms a weak ascent together with the first element of the next block $B_{t+1}$). Then
$$\TPS(x)=\mathcal{R} (B_1) \cdots \mathcal{R}(B_k).$$
\end{lemma}

\begin{theorem}\label{tortoise_pop_char}
Let $\pi$ a Cayley permutation. Then $\pi$ is tortoise-popstack sortable if and only if $\pi$ avoids $231$, $312$, $221$ and $211$.
\end{theorem}
\proof The proof is similar to that of Theorem~\ref{hare_pop_char}. It is not difficult to show that if $\pi$ contains an occurrence of either $231$, $312$, $221$ or $211$, then $\pi$ is not tortoise pop-stack sortable. We leave the details to the reader.\\
Conversely, suppose that $\pi$ is not tortoise-popstack sortable. Then $\TPS(\pi)=\mathcal{R} (B_1) \cdots \mathcal{R}(B_k)$ and $\TPS(\pi)$ contains at least one strict descent $a>b$. Suppose $a,b$ is the leftmost strict descent in $\TPS(\pi)$. Again it has to be $a \in B_i$ and $b \in B_{i+1}$ for some $i$, due to Lemma~\ref{blocks_dec_tortoise}. Denote with $u$ the last element of $B_i$ and with $v$ the first element of $B_{i+1}$. Here it must be $a > u$, $v > b$ and $u \le v$ by Lemma~\ref{blocks_dec_tortoise} (see Figure~\ref{blocks_fig}). Suppose that $B_{i}$ is a singleton and thus $a=u$. Then $v \neq b$, since $v \ge a$, whereas $b<a$. Now, if $a<v$, then $avb$ is an occurrence of $231$. Otherwise, if $a=v$, then $avb$ is an occurrence of $221$. Otherwise, suppose that $B_{i}$ is not a singleton and thus $a \neq u$. If $a>v$, then $auv$ is an occurrence of either $312$, if $u<v$, or $211$, if $u=v$. If $a<v$, then $v \neq b$, since $v>a$ and $b<a$, and $avb$ is an occurrence of $231$. Finally, suppose that $a=v$. Note that again $v \neq b$, since $v=a$ and $b<a$. Therefore, if $b=u$, then $aub$ is an occurrence of $211$. If $b<u$, then $uvb$ is an occurrence of either $231$, if $u<v$, or $221$, if $u=v$. Instead, if $b>u$, then $aub$ is an occurrence of $312$.
\cvd

Let us now enumerate tortoise pop-stack sortable Cayley permutations. First a geometrical description. Write again $\pi=B_1 B_2 \cdots B_k$, where each block $B_i$ is maximally strictly decreasing. Denote by $m_i$ the first element of $B_i$ and let $B_i=m_i A_i$, where $A_i$ contains the remaining elements of $B_i$. Suppose that $\pi$ is tortoise pop-stack sortable. Then:

\begin{enumerate}
\item $m_i \le m_{i+1}$ for each $i$. Otherwise, suppose by contradiction that $m_i>m_{i+1}$. Let $x$ the last element of $B_i$. Then $x \le m_{i+1}$, therefore $x \neq m_i$ and $m_i x m_{i+1}$ is an occurrence of $312$, against Theorem~\ref{tortoise_pop_char}.

\item $A_i<A_{i+1}$ for each $i$. In other words, $x<y$ for each $x \in A_i$ and $y \in A_{i+1}$. Otherwise, if $x=y$, then $m_i x y$ is an occurrence of $211$, against Theorem~\ref{tortoise_pop_char}. Instead, if $x>y$, then $x m_{i+1} y$ is an occurrence of $231$, which is impossible due to the same result.

\item $y \ge m_i$ for each $y \in A_{i+1}$. Otherwise $m_i m_{i+1} y$ is an occurrence of either $231$, if $m_i<m_{i+1}$, or $221$, if $m_i=m_{i+1}$. In both cases this is impossible due to Theorem~\ref{tortoise_pop_char}.
\end{enumerate}

Denote by $f_{n,k}$ the number of tortoise pop-stack sortable Cayley permutations of length $n$ and with $k$ maximally strictly decreasing blocks. As a consequence of what said above, each of these Cayley permutations is determined uniquely by choosing:

\begin{itemize}
\item the length of each block, which can be done in $\binom{n-1}{k-1}$ distinct ways, and
\item whether the first element of a block is equal to or greater than the first element of the previous block. Equivalently, whether $m_{i+1}=m_i$ or $m_{i+1}=m_i+1$, for each $i \ge 2$. 
\end{itemize}

Therefore $f_{n,k} = \binom{n-1}{k-1} 2^{k-1}$.

\begin{corollary}
For each $n \ge 1$, there are $3^{n-1}$ tortoise pop-stack sortable permutations of length $n$.
\end{corollary}
\proof Let $f_n$ be the number of tortoise pop-stack sortable permutations of length $n$. We have
\begin{equation}
\begin{split}
3^{n-1}=(2+1)^{n-1}=\sum_{j=0}^{n-1} \binom{n-1}{j} 2^j &=\\
\sum_{k=1}^n \binom{n-1}{k-1} 2^{k-1} = \sum_{k=1}^n f_{n,k} =f_n &\\
\end{split}
\end{equation}
\cvd

The enumeration of hare pop-stack sortable Cayley permutations, or equivalently the set of Cayley permutations avoiding $231$, $312$ and $2121$, is rather more complicated, thus we leave it for a future work.

\begin{openpr}
Enumerate the hare pop-stack sortable Cayley permutations. The sequence starts $1,3,11,41,151,553$ and it does not match any sequence in the OEIS~\cite{Sl}.
\end{openpr}

%

\subsubsection*{Acknowledgements}
The author would like to sincerely thank Anders Claesson for the fruitful discussions and the precious data he provided.
%

\end{document}